\documentclass[11pt]{elsarticle}

\usepackage[utf8]{inputenc}
\usepackage[english]{babel}
\usepackage{enumerate}
\usepackage{latexsym}

\usepackage{multicol}
\usepackage{tikz}
\usepackage{float} 

\usepackage{amsthm}
\usepackage{amsmath}
\usepackage{amssymb}
\usepackage{aliascnt} 

\usepackage[bookmarks]{hyperref} 
\hypersetup{
	pdftitle={2-switch-degree Classification of Split Graphs},
	pdfauthor={V.N. Schvollner},
	colorlinks=true,
	linkcolor=blue,
	citecolor=red,
	filecolor=cyan,
	urlcolor=magenta
}

\usetikzlibrary{babel,decorations.pathreplacing,calc,positioning,shapes.geometric}

\usepackage{cleveref}

\definecolor{10}{RGB}{115,59,171}
\definecolor{8}{RGB}{212,122,240}
\definecolor{7}{RGB}{99,212,119}
\definecolor{6}{RGB}{183,240,164}
\definecolor{D}{RGB}{255,162,79}
\definecolor{E}{RGB}{255,84,0}
\definecolor{F}{RGB}{158,248,255}
\definecolor{G}{RGB}{128,135,255}
\definecolor{I}{RGB}{187,255,0}
\definecolor{A}{cmyk}{.9,.05,.4,0}
\definecolor{B}{RGB}{150,30,150}
\definecolor{C}{RGB}{186,155,189}
\definecolor{9}{RGB}{0,180,60}
\definecolor{0}{RGB}{30,123,191}
\definecolor{1}{RGB}{255,113,102}
\definecolor{2}{RGB}{41,199,92}
\definecolor{3}{RGB}{242,207,16}
\definecolor{5}{RGB}{255,15,154}
\definecolor{4}{rgb}{.8,0,.8}

\definecolor{Red}{rgb}{1,0.4,0.4}
\definecolor{Green}{rgb}{.1,.5,.1}
\definecolor{Blue}{rgb}{.1,.1,.5}
\definecolor{blue}{RGB}{0,0,255}
\definecolor{Yellow}{rgb}{.8,.4,0}
\definecolor{X}{rgb}{.8,.4,0}
\definecolor{H}{rgb}{0,0,1}
\definecolor{light}{rgb}{.67,.84,.90}
\definecolor{Cyan}{rgb}{0,1,1}
\definecolor{Purple}{rgb}{.5,0,.5}
\definecolor{Purple2}{rgb}{.5,.2,.5}
\definecolor{white}{rgb}{1.0,1.0,1.0}
\definecolor{Purple2}{rgb}{.8,.4,0}
\definecolor{Amarillo}{RGB}{225,191,73}
\definecolor{Celeste}{RGB}{117,170,219}
\definecolor{Castano}{RGB}{232,53,17}
\definecolor{Black}{RGB}{0,0,0}
\definecolor{White}{RGB}{255,255,255}
\definecolor{gris}{rgb}{.5,.5,.5}

\newtheorem{theorem}{Theorem}[section]

\newaliascnt{corollary}{theorem}
\newtheorem{corollary}[corollary]{Corollary}
\aliascntresetthe{corollary}

\newaliascnt{lemma}{theorem}
\newtheorem{lemma}[lemma]{Lemma}
\aliascntresetthe{lemma}

\newaliascnt{proposition}{theorem}
\newtheorem{proposition}[proposition]{Proposition}
\aliascntresetthe{proposition}

\DeclareMathOperator{\act}{act}

\pgfdeclarelayer{background2}
\pgfdeclarelayer{background}
\pgfdeclarelayer{foreground}
\pgfsetlayers{background2,background,main,foreground}

\begin{document}
	\begin{frontmatter}
		\title {
			Simple factor graphs associated with split graphs
		}
		\address[IMASL]{Instituto de Matem\'atica Aplicada San Luis, Universidad Nacional de San Luis and CONICET, San Luis, Argentina.}
		\address[DEPTO]{Departamento de Matem\'atica, Universidad Nacional de San Luis, San Luis, Argentina.}

		\author[IMASL]{Victor N. Schvöllner}\ead{vnsi9m6@gmail.com} 
		\author[IMASL,DEPTO]{Adri\'an Pastine}\ead{agpastine@unsl.edu.ar}

		\begin{abstract}
			We introduce and study a loopless multigraph associated with a split graph $(S,K,I)$: its factor graph $\Phi(S,K,I)$, whose vertex set is $I$ and whose edges, counted with multiplicity, record the 2-switches acting on $S$. In contrast with the $A_4$-structure $\mathcal{A}_4(S)$ introduced by Barrus and West, the factor graph is supported on the independent set alone; it is therefore more compact and, although it retains less information, it still determines the 2-switch-degree of $S$. We show that an active split graph is indecomposable if and only if its factor graph is connected, an analogue for split graphs of the Barrus--West characterization of indecomposable graphs. Finally, we describe completely the active split graphs whose factor graph is simple and connected.
		\end{abstract}
	\end{frontmatter}	
	
	
\section{Introduction}

For specific graph theoretical concepts not defined here, we refer the reader to \cite{chartrand2010graphs}. Let $G$ be a graph. We use the notation $V(G)$ and $E(G)$ to refer to the vertex set and the edge set of $G$, respectively. The complement of $G$ is denoted by $\overline{G}$. On the other hand, $H\prec G$ means that $H$ is an induced subgraph of $G$. If $V'=\{v_1,\ldots,v_n\}\subseteq V(G)$, we denote by $G[V']$ or $G[v_1,\ldots,v_n]$ the subgraph of $G$ induced by the vertices in $V'$. We denote a path by the juxtaposition of its vertices in order; for instance, $uxyv$ is the $P_4$ with edges $ux,xy,yv$ (and non-edges $uy,uv,xv$ when the $P_4$ is induced). The set $N_{G}(v)=\{x\in V(G): vx\in E(G)\}$ is the \emph{neighborhood} of $v$ in $G$. The \emph{degree} of $v$ in $G$ (i.e., $|N_G(v)|$) is denoted by $\deg_G(v)$. We usually write $N_v$ and $d_v$, instead of $N_G(v)$ and $\deg_G(v)$ respectively, when $G$ is clear from the context. The \emph{degree sequence} of $G$ is the $|V(G)|$-tuple whose $v$-coordinate is $\deg_G(v)$. The \emph{clique number} of $G$, denoted by $\omega(G)$, is the size of a maximum clique in $G$. The \emph{independence number} of $G$, denoted by $\alpha(G)$, is the size of a maximum independent set of $G$. Clearly, $\omega(G)=\alpha(\overline{G})$, for every graph $G$. The notation $G\approx H$ means that the graphs $G$ and $H$ are isomorphic. For $n\in\mathbb{N}=\{1,2,\ldots\}$, $[n]$ means $\{1,\ldots,n\}$.\\

Let $G$ be a graph containing four distinct vertices $a,b,c,d$ such that $ab,cd \in E(G)$ and $ac,bd \notin E(G)$. A \emph{2-switch} on $G$ is an operation that removes the edges $ab$ and $cd$ from $G$ and adds the edges $ac$ and $bd$ (see \cite{chartrand2010graphs}, page 23). Directly from the definition, it follows that 2-switch preserves the degree sequence. Moreover, it replaces edges between vertices $a,b,c,d$ in $G$ if and only if $G[a,b,c,d]\approx P_4,C_4$ or $2K_2$. The \emph{realization graph} is the graph whose vertex set is the family of all labeled graphs with the same degree sequence, where two vertices (graphs) $G,H$ are adjacent if and only if $H=\tau(G)$ for some 2-switch $\tau$ (see \cite{pastine20252}, pages 3 and 4). Given any two graphs $G,H$ with the same degree sequence, a classical theorem of Fulkerson, Hoffman and McAndrew (see \cite{fulkerson1965}) states that $G$ can be transformed into $H$ by a sequence of 2-switches (see also \cite{chartrand2010graphs}, page 24). As a consequence, every realization graph is connected.

Two fundamental results illustrate the importance of 2-switch transformations and motivate this article: the canonical decomposition theorem by R. Tyshkevich \cite{tyshkevich2000decomposition}, and the characterization by Barrus and West \cite{barrus.west.A4} of indecomposable graphs. In order to properly understand the content of both, we need to clarify these preliminary concepts: split graph, Tyshkevich composition, and indecomposable graph. 

A graph $G$ is said to be \emph{split} if $V(G)\neq\varnothing$ is the disjoint union of two sets, $V(G)=K\dot{\cup} I$, where $K$ is a clique and $I$ is an independent set. In such a case, we say that the pair $(K,I)$ is a \emph{bipartition} for $G$. The notation $(S,K,I)$ means that we are considering the split graph $S$ together with a bipartition $(K,I)$ for it. Further information about the structure of split graphs can be found in \cite{splitnordhausgaddum,hammer1977split}. If $(S,K,I)$ is a split graph and $G$ is a graph disjoint from $S$, the \emph{Tyshkevich composition} $(S,K,I)\circ G$ of $(S,K,I)$ and $G$ is defined as the graph whose vertex set is $V(S)\cup V(G)$ and whose edge set is $E(S)\cup E(G)\cup\{xy:x\in K,y\in V(G)\}$. This operation was introduced by R. Tyshkevich in \cite{tyshkevich2000decomposition}. In general, $S\circ G$ is a split graph if and only if $G$ is. A graph $G$ is \emph{decomposable} if there exist graphs $S, H$ such that $S\circ H=G$ and $V(S), V(H)\neq\varnothing$. Otherwise, $G$ is said to be \emph{indecomposable}.

\begin{theorem}[Tyshkevich, \cite{tyshkevich2000decomposition}]
	\label{tyshk.decomp}
	Every graph $G$ can be written as a (Tyshkevich) composition of indecomposable graphs: $G=G_n \circ \ldots\circ G_2 \circ G_1$,	where $G_2,\ldots,G_n$ are split. When $|V(G_r)|\geq 1$ for all $r$, this decomposition is unique up to isomorphism, and neither the order of the factors nor the choice of bipartitions of the split factors can vary.
\end{theorem}

Notice that \Cref{tyshk.decomp} highlights the importance of indecomposable split graphs as fundamental building blocks. 

Following \cite{barrus.west.A4}, four distinct vertices $a,b,c,d$ of $G$ form an \emph{alternating $4$-cycle} (or \emph{$A_4$-cycle}) if $G[a,b,c,d]\approx P_4,C_4$ or $2K_2$. Thus the $A_4$-cycles of $G$ are exactly the $4$-sets on which a $2$-switch can be performed. Barrus and West studied the \emph{$A_4$-structure} of $G$: the $4$-uniform hypergraph $\mathcal{A}_4(G)$ with $V(\mathcal{A}_4(G))=V(G)$ whose edges are the $A_4$-cycles of $G$. They proved that $\mathcal{A}_4(G)$ decomposes as the disjoint union of the $A_4$-structures of the Tyshkevich factors of $G$; as a consequence, $G$ is indecomposable if and only if $\mathcal{A}_4(G)$ is connected (see Theorem 3.2 of \cite{barrus.west.A4}, page 7). The \emph{$2$-section} of $\mathcal{A}_4(G)$ is the graph $A_4(G)$ with $V(A_4(G))=V(G)$ in which $uv\in E(A_4(G))$ if and only if some $A_4$-cycle of $G$ contains both $u$ and $v$, i.e. if and only if some $2$-switch on $G$ involves both $u$ and $v$. A hypergraph is connected if and only if its $2$-section is connected, because paths in the 2-section correspond exactly to chains of consecutively intersecting hyperedges. Hence, we have the following result:

\begin{corollary}
	\label{indecomp.characterization}
	Let $G$ be a graph. If $G_n\circ\ldots\circ G_1$ is the Tyshkevich decomposition of $G$, then
	\[ A_4(G)=\dot{\bigcup}_{i=1}^n A_4(G_i). \]
	In particular, $G$ is indecomposable if and only if $A_4(G)$ is connected.
\end{corollary}

\begin{proof}
	It follows from Theorem 3.2 of \cite{barrus.west.A4} and the previous discussion.
\end{proof}

The passage from $\mathcal{A}_4(G)$ to its $2$-section $A_4(G)$ is lossy: each $A_4$-cycle of $G$ collapses to a $4$-clique of $A_4(G)$, so $A_4(G)$ retains which pairs of vertices participate together in some $2$-switch but no longer records how many, nor the identity of the $A_4$-cycles themselves. Furthermore, a $4$-clique of $A_4(G)$ need not come from an $A_4$-cycle. To see that, let $(S,K,I)$ be the split graph with $I=[4]$, $K=\{5,6,7,8\}$ and $N_v=\{v+4\}$ for $v\in I$. For distinct $x,y\in K$, we have $S[x-4,x,y,y-4]\approx P_4$, so $A_4(S)[K]\approx K_4$; yet $S[K]\approx K_4$ is not an $A_4$-cycle. This becomes an issue as soon as we want to count $2$-switches acting on $S$, i.e. the \emph{$2$-switch-degree} of $S$.

Computing the $2$-switch-degree (or simply, the \emph{degree}) of a graph $G$, denoted by $\deg(G)$, is a general problem first introduced and studied by Schv\"ollner and Pastine in \cite{pastine20252}. In that work, the authors recognize, among other things, the importance of $\deg(G)$ in understanding the structure of the realization graph associated with the degree sequence of $G$. Furthermore, for a split graph $S$, $\deg(S)$ equals the number of induced $P_4$'s of $S$ (see \cite{pastine20252}, page~15), because split graphs are precisely those that do not contain induced $C_4,C_5$ or $2K_2$ (see \cite{hammer1977split}). Based on all these observations and inspired by Barrus and West's $A_4$-structure, in this article we introduce the concept of factor graph associated with a split graph.  

The \emph{factor graph} of $(S,K,I)$ is the loopless multigraph $\Phi=\Phi(S,K,I)$ with $V(\Phi)=I$, such that there is an edge between $u$ and $v$ for each $P_4\prec S$ containing both $u$ and $v$. Thus, $E(\Phi)$ is a multiset, and so every edge $uv$ of $\Phi$ has a multiplicity, that we denote by $\sigma_{S}(uv)$. As usual, we can simply write $\sigma_{uv}$ if $S$ is clear from the context. If $\sigma_{uv}\in\{0,1\}$ for all $\{u,v\}\subseteq I$, we say that $\Phi$ is \emph{simple}. A concrete example of a split graph $S$ with its associated $A_4(S)$ and $\Phi$ is shown in \Cref{fig:ejemplo_A4_Phi} (where the $K$-vertices are filled with black).
\begin{figure}[htbp]
	\centering
	\begin{tikzpicture}[
		scale=1,
		K/.style={circle,draw,fill=black,minimum size=5.5pt,inner sep=0pt},
		I/.style={circle,draw,fill=white,minimum size=5.5pt,inner sep=0pt},
		every edge/.style={draw,line width=.5pt}
		]
		
		\begin{scope}
			\node[K] (a) at (0,0)      {}; \node[below left]  at (a) {$4$};
			\node[K] (b) at (0.7,1.15)  {}; \node[right=2pt]       at (b) {$5$};
			\node[K] (c) at (1.55,0.1) {}; \node[below right]        at (c) {$6$};
			\node[I] (u1) at (-1.55,0.2){}; \node[below=2pt]        at (u1){$1$};
			\node[I] (u2) at (-0.9,1.2) {}; \node[left=2pt]  at (u2){$2$};
			\node[I] (u3) at (0.7,-1.15){}; \node[below right]       at (u3){$3$};
			\path (a) edge (b) (b) edge (c) (a) edge (c);
			\path (u1) edge (b);
			\path (u2) edge (a) (u2) edge (b);
			\path (u3) edge (a) (u3) edge (c);
			\node at (0.1,-2) {$S$};
		\end{scope}
		
		\begin{scope}[xshift=4.6cm,yshift=-0.05cm]
			\node[I] (A1) at (-0.65,1.13) {}; \node[above left]  at (A1){$1$};
			\node[I] (A2) at ( 0.65,1.13) {}; \node[above right] at (A2){$2$};
			\node[K] (Ab) at ( 1.30,0)    {}; \node[right=2pt]       at (Ab){$5$};
			\node[K] (Aa) at ( 0.65,-1.13){}; \node[below right] at (Aa){$4$};
			\node[I] (A3) at (-0.65,-1.13){}; \node[below left]  at (A3){$3$};
			\node[K] (Ac) at (-1.30,0)    {}; \node[left=2pt]        at (Ac){$6$};
			\path (A1) edge (A3) (A1) edge (Aa) (A1) edge (Ab) (A1) edge (Ac);
			\path (A2) edge (A3) (A2) edge (Ab) (A2) edge (Ac);
			\path (A3) edge (Aa) (A3) edge (Ab) (A3) edge (Ac);
			\path (Aa) edge (Ab) (Ab) edge (Ac);
			\node at (0,-2) {$A_4(S)$};
		\end{scope}
		
		\begin{scope}[xshift=8.0cm,yshift=-0.75cm]
			\node[I] (P1) at (0,1.15)   {}; \node[above left]  at (P1){$1$};
			\node[I] (P2) at (1.5,1.15) {}; \node[above left] at (P2){$2$};
			\node[I] (P3) at (0,-0.35)  {}; \node[below left]  at (P3){$3$};
			\path (P1) edge[bend left=18]  (P3);
			\path (P1) edge[bend right=18] (P3);
			\path (P2) edge (P3);
			\node at (0.75,-1.3) {$\Phi(S,K,I)$};
		\end{scope}
		
	\end{tikzpicture}
	\caption{A split graph $(S,K,I)$, where $I=[3]$, together with $A_4(S)$ and $\Phi(S,K,I)$.}
	\label{fig:ejemplo_A4_Phi}
\end{figure}

By definition, each edge of $\Phi$ records an induced $P_4$ of $S$ whose two endpoints are $u,v\in I$. Recall that split graphs have no induced $C_4$ nor $2K_2$, so every $A_4$-cycle of $S$ induces a $P_4$. Moreover, if $uxyv\prec S$, then $u,v\in I$ and $x,y\in K$.
Thus the ``endpoint map'' $P_4\mapsto\{u,v\}$ is a bijection between the induced $P_4$'s of $S$ and the edges of $\Phi$, counted with multiplicity. As each such $P_4$ carries a single $2$-switch, we obtain, for every fixed bipartition $(K,I)$,
\[ \deg(S)=|E(\mathcal{A}_4(S))|=|E(\Phi)|=\sum_{\{u,v\}\subseteq I}\sigma_S(uv). \]
What is not immediate is that this edge-multiset is the same for every bipartition. This is settled in \Cref{sec:Phi.y.biparticiones}, where we show that the dependence of $\Phi$ on $(K,I)$ is confined to isolated vertices, namely those members of $I$ that lie in no induced $P_4$ of $S$. Consequently, all factor graphs of $S$ share a common ``core'', obtained by discarding their isolated vertices and any two of them differ by at most one isolated vertex (\Cref{Phi.bipartition.change,Phi.dos.niveles}). In particular, $E(\Phi)$ is an invariant of $S$, and $\deg(S)=|E(\Phi)|$ irrespective of the bipartition. Accordingly, from now on we may write $\Phi(S)$ instead of $\Phi(S,K,I)$ when $S$ has a unique bipartition.

In passing from $\mathcal{A}_4(S)$ to $\Phi$ we forget which pair of $K$-vertices each induced $P_4$ contains, and $\Phi$ cannot recover how many induced $P_4$'s of $S$ contain a given $x\in K$, or a given pair $x,y\in K$. However, what it does retain is exactly the data we need in this article:
\begin{enumerate}
	\item the total number of $2$-switches acting on $S$ (i.e., $\deg(S)$), which is $|E(\Phi)|$ counting multiplicities;
	
	\item the number of $2$-switches involving a prescribed pair $\{u,v\}\subseteq I$, which is $\sigma_S(uv)$;
	
	\item the number of $2$-switches involving a prescribed vertex $v\in I$, which is $\deg_{\Phi}(v)$ counting multiplicities;
	
	\item whether $S$ is indecomposable.
\end{enumerate}
Compared with the $2$-section $A_4(S)$, the gain is clear: neither $\deg(S)$ nor $\sigma_{uv}$ is recoverable from $A_4(S)$ alone; $\Phi$ records both. Moreover, ignoring multiplicities, $\Phi \prec A_4(S)$.

In \cite{pastine20252} (pages 6, 8 and 12) Pastine and Schvöllner introduce the concepts of active vertex, active graph and prime graph. We say that a vertex $v$ of a graph $G$ is \emph{active} in $G$ if $v$ participates in some 2-switch over $G$. Otherwise, $v$ is \emph{inactive} in $G$. The set of all active vertices of $G$ is denoted by $\act(G)$. $G$ is said to be \emph{active} if $\act(G)=V(G)$. When $G$ is active and indecomposable, we say that $G$ is \emph{prime}. As noticed in \cite{pastine20252} (page 12), prime graphs are exactly indecomposable nontrivial graphs. 

\Cref{sec:relacion.Phi-A_4} relates the connectedness of $\Phi$ to that of $A_4(S)$. The key observation is that, ignoring multiplicities, $\Phi$ is the incomparability graph of $\{N_S(v):v\in I\}$; from it we deduce that a connected $A_4(S)$ forces $\Phi$ to be connected. As a consequence, an active split graph $S$ is prime if and only if $\Phi(S)$ is connected (\Cref{S.primo.iff.Phi(S).conexo}), the analogue of \Cref{indecomp.characterization} for split graphs.

\Cref{sec:homogen_split} studies how $\sigma_{uv}$ reflects the inclusion and intersection of $N_u$ and $N_v$, and introduces the concept of $\delta$-homogeneous split graphs, whose $I$-vertices all have degree $\delta$. When $\Phi$ is simple and connected, $S$ turns out to be $\delta$-homogeneous and balanced (\Cref{Phi.simple.conexo.1}), which motivates our focus on the simple case in \Cref{sec:simple_complete_Phi,sec:simple_Phi}.

In \Cref{sec:simple_complete_Phi}, we study the multigraph $\Phi$ under the assumption that it is simple and complete. In such a case, using a result from \cite{automorphism.johnson.graph} about a certain intersecting family of sets, we describe the structure of $S$, and we show that for active $S$ this occurs exactly when $|K|=|I|$ and the edges between $K$ and $I$ form a perfect matching or the complement of one (\Cref{Phi.simple.completo.implica...,act.Phi.simple.compl.iff...}).

In \Cref{sec:simple_Phi} we assume $S$ active and $\Phi$ simple and connected, and we describe every such $S$ completely: each $v\in I$ has a unique neighbor in $K$, or each $v\in I$ has a unique non-neighbor in $K$, and $\Phi$ is the complete multipartite graph whose parts are the classes of twins of $I$ (\Cref{caract.Phi.simple.conexo}); the simple and complete case of \Cref{sec:simple_complete_Phi} is precisely the twin-free specialization. From this description we read off $\omega(S)=\omega(\Phi)=|K|\leq|I|=\alpha(S)$, with equality exactly when $\Phi$ is complete, together with an exact formula for $\deg(S)$ (\Cref{Phi.simple.conexo.deg.cotas}). We close by showing that $\Phi(\overline{S})$ is always complete, and that it is simple exactly when $\Phi(S)$ is (\Cref{Phi.complement.not.simple}).

The main takeaways of this article are two: the characterization of primality for active split graphs in terms of the connectedness of $\Phi$, and the complete description of the active split graphs whose factor graph is simple and connected. The remaining results are the tools developed to reach them.


\section{Factor graphs and bipartitions}
\label{sec:Phi.y.biparticiones}

The multigraph $\Phi$ is built from a bipartition for $S$, not from $S$ alone, so before developing its theory we must clarify to what extent it depends on that choice. We show that the dependence is confined to isolated vertices.

A split graph $(S,K,I)$ is \emph{balanced} if $|K|=\omega(S)$ and $|I|=\alpha(S)$, and \emph{unbalanced} otherwise. A vertex $w$ of $S$ is \emph{swing} if there exist bipartitions $(K,I)$ and $(K',I')$ for $S$ with $w\in K$ and $w\in I'$; we write $W(S)$ for the set of swing vertices, an invariant of $S$. The term \emph{swing} first appeared in \cite{whitman2020split}, with a bipartition-dependent definition later refined to the present one in Section~4 of \cite{jaume2025nullspace}. From there, we highlight that $S$ is balanced if and only if $W(S)=\varnothing$ if and only if $S$ has a unique bipartition. Hence, the symbol $\Phi(S)$ is unambiguous when $S$ is balanced. In particular, the same holds if $S$ is active, because every active split graph is balanced, by Proposition 4.1(7) of \cite{pastine20252}. 

A vertex is inactive if and only if it lies in no induced $P_4$ of $S$, so whenever such a vertex belongs to $I$, it is isolated in $\Phi(S,K,I)$. In particular, the members of $W(S)\cap I$ are isolated in $\Phi(S,K,I)$, because swing vertices are inactive, by Proposition 4.1(6) of \cite{pastine20252}. Following Section~5 of \cite{pastine20252}, the \emph{active part} of $S$ is the subgraph $S^*=S[\act(S)]$. Since $S^*$ is itself active, $\Phi(S^*)$ is unambiguous. When $S$ is unbalanced, it admits more than one bipartition, but the next proposition shows that all the resulting factor graphs agree once their isolated vertices are removed, sharing the subgraph $\Phi(S^*)$.

\begin{proposition}\label{Phi.bipartition.change}
	Let $S$ be a split graph and let $A$ be any set with $W(S)\subseteq A\subseteq V(S)-\act(S)$. Then, for every bipartition	$(K,I)$ of $S$, the multigraph $\Phi(S,K,I)-(A\cap I)$ is the same. Furthermore, each deleted vertex is isolated in $\Phi(S,K,I)$. In particular, any two factor graphs of $S$ have the same edges with the same multiplicities. Taking $A=V(S)-\act(S)$ yields
	\[
	\Phi(S,K,I)-\bigl(I-\act(S)\bigr)=\Phi(S^*),
	\]
	so all factor graphs of $S$ share the subgraph $\Phi(S^*)$, independently of the bipartition for $S$.
\end{proposition}

\begin{proof}
	By the discussion above, we know that $A\cap I$ consists of vertices that are isolated in $\Phi(S,K,I)$, and deleting it preserves all edges and multiplicities. The residual vertex set is $I-A$. If $(K',I')$ is another bipartition, then $I\bigtriangleup I'\subseteq W(S)\subseteq A$, so $I- A=I'- A$; and for $u,v$ in this common set, $\sigma_{uv}$ is the number of induced $P_4$'s of $S$ with endpoints $u$ and $v$, a quantity independent of the bipartition. Hence $\Phi(S,K,I)-(A\cap I)$ does not depend on $(K,I)$. Finally, take $A=V(S)-\act(S)$: the residual vertices are $\act(S)\cap I$, and every induced $P_4$ of $S$ has all four of its vertices active, hence it is induced in $S^*$. Therefore $\sigma_S(uv)=\sigma_{S^*}(uv)$ and the residual multigraph is $\Phi(S^*)$.
\end{proof}

For an unbalanced $S$, the exact number of isolated vertices depends only on whether $I$ is a maximum independent set of $S$.

\begin{proposition}\label{Phi.dos.niveles}
	Let $(S,K,I)$ be an unbalanced split graph and set $p=|\act(S)\cap I|$. Then,	$|I|\in\{\alpha(S),\alpha(S)-1\}$ and:
	\begin{enumerate}
		\item if $|I|=\alpha(S)-1$, then $\Phi(S,K,I)=\Phi(S^*)\;\dot\cup\;\overline{K_{\,\alpha(S)-1-p}}\, ;$
		
		\item if $|I|=\alpha(S)$, then $\Phi(S,K,I)=\Phi(S^*)\;\dot\cup\;\overline{K_{\,\alpha(S)-p}}\, .$
	\end{enumerate}
	In both cases the vertex set of the edgeless part is $I-\act(S)$. All bipartitions whose independent part has the same size yield isomorphic factor graphs, and a	factor graph of case (2) has exactly one more isolated vertex than a factor graph of case (1).
\end{proposition}
\begin{proof}
	It is a well-known fact that every bipartition $(K,I)$ for a split graph $S$ satisfies $\alpha(S)-1\le|I|\le\alpha(S)$ (see \cite{splitnordhausgaddum} or \cite{hammer1977split}). By \Cref{Phi.bipartition.change}, a vertex of $I$ is isolated in $\Phi(S,K,I)$ if and only if it is inactive, and deleting the isolated vertices leaves $\Phi(S^*)$; hence $\Phi(S,K,I)=\Phi(S^*)\;\dot\cup\;\overline{K_q}$ with $q=|I|-|\act(S)\cap I|$. An active vertex is never swing, so it lies on the same side of every bipartition; thus $\act(S)\cap I$ is independent of $(K,I)$, and $q=|I|-p$ equals $\alpha(S)-1-p$ in case~(1) and $\alpha(S)-p$ in case~(2).
\end{proof}

In view of \Cref{Phi.bipartition.change,Phi.dos.niveles}, the choice of bipartition affects $\Phi$ only through isolated inactive vertices. As an example, let $S$ be the split graph with $V(S)=[5]$ and $E(S)=\{34,35,45,13,24\}$. Its bipartitions are $(K,I)=(\{3,4,5\},[2])$ and $(K',I')=(\{3,4\},\{1,2,5\})$; $W(S)=\{5\}$, so $S$ is unbalanced with $\alpha(S)=3$ and $p=2$. Therefore, $\Phi(S,K,I)=\Phi(S^*)\approx K_2$ and $\Phi(S,K',I')=\Phi(S^*)\dot\cup K_1$. The two factor graphs differ exactly by the isolated swing vertex $5$, as predicted by \Cref{Phi.dos.niveles}.

\begin{corollary}
	Let $S$ be a split graph and let $(K,I)$ be any bipartition for $S$. If $\Phi=\Phi(S,K,I)$, then $E(\Phi)=E(\Phi(S^*))$ and $|E(\Phi)|=\deg(S)$. In other words, the edge-multiset of $\Phi$ does not depend on the bipartition. 
\end{corollary}

\begin{proof}
	It follows immediately from \Cref{Phi.dos.niveles}.
\end{proof}


\section{Relationship between $\Phi$ and $A_4(S)$} \label{sec:relacion.Phi-A_4}


If $A, B$ are sets, we say that $A$ and $B$ are \emph{comparable} if $A\subseteq B$ or $B\subseteq A$, and \emph{incomparable} otherwise. We denote $|N_S(u)\cap N_S(v)|$ by $\eta_{S}(uv)$. As usual, if no ambiguity arises in the context, we may write simply $\eta_{uv}$.

\begin{proposition}
	\label{prop.incomparabilidad}
	Let $(S,K,I)$ be a split graph and let $u,v\in I$ be distinct.
	\begin{enumerate}
		\item $\sigma_{uv}=|N_u-N_v|\,|N_v-N_u|=(d_u-\eta_{uv})(d_v-\eta_{uv})$;
		
		\item $uv\in E(\Phi(S,K,I))$ if and only if $N_u$ and $N_v$ are incomparable;
    \end{enumerate}
\end{proposition}

\begin{proof}
	\begin{enumerate}[(1).]
		\item By definition, $\sigma_S(uv)=|\{H\prec S: u,v\in V(H), H\approx P_4\}|$. Since $u,v\in I$, such an $H$ has $u$ and $v$ as its endpoints and its two internal vertices in $K$; writing them as $x$ (adjacent to $u$) and $y$ (adjacent to $v$), the non-adjacencies $uy,xv\notin E(S)$ say exactly that $x\in N_u-N_v$ and $y\in N_v-N_u$. Conversely, given any $x\in N_u-N_v$ and any $y\in N_v-N_u$, we have $x\neq y$, $xy\in E(S)$ because $x,y\in K$, and $uv\notin E(S)$ because $u,v\in I$; hence $S[u,x,y,v]\approx P_4$. The assignment $H\mapsto (x,y)$ is therefore a bijection onto $(N_u-N_v)\times(N_v-N_u)$, which gives the first equality. The second follows from $|N_u-N_v|=d_u-\eta_{uv}$ and $|N_v-N_u|=d_v-\eta_{uv}$.
		
		\item By (1), $\sigma_S(uv)\neq 0$ if and only if $N_u-N_v\neq\varnothing$ and $N_v-N_u\neq\varnothing$. \qedhere
	\end{enumerate}
\end{proof}

In \Cref{ejemplo.cantidad.P4.ind.split} we see an example of an active split graph $S$ where $\deg(S)$ is computed by using \Cref{prop.incomparabilidad}(1). 
\begin{figure}[H]
	\centering
	\begin{tikzpicture}[
		scale=1,
		K/.style={circle,draw,fill=black,minimum size=5.5pt,inner sep=0pt},
		I/.style={circle,draw,fill=white,minimum size=5.5pt,inner sep=0pt},
		every edge/.style={draw,line width=.5pt},
		Fedge/.style={draw,line width=1.8pt,color=black!45}
		]
		
		\begin{scope}
			\node[K] (n1) at (0,1.2)        {}; \node[above] at (n1){};
			\node[K] (n2) at (1.141,0.371)  {}; \node[above]  at (n2){};
			\node[K] (n3) at (0.705,-0.971) {}; \node[below right] at (n3){};
			\node[K] (n4) at (-0.705,-0.971){}; \node[below left] at (n4){};
			\node[K] (n5) at (-1.141,0.371) {}; \node[below left] at (n5){};
			\node[I] (a) at (-2.6,0.45) {}; \node[above=2pt] at (a){$1$};
			\node[I] (b) at ( 2.6,0.55) {}; \node[above right]at (b){$2$};
			\node[I] (c) at ( 2.6,-0.5) {}; \node[below right]at (c){$3$};
			\path (n1) edge (n2) (n2) edge (n3) (n3) edge (n4) (n4) edge (n5) (n5) edge (n1);
			\path (n1) edge (n3) (n1) edge (n4) (n2) edge (n4) (n2) edge (n5) (n3) edge (n5);
			\path (a) edge (n1) (a) edge (n4) (a) edge (n5);
			\path (b) edge (n1) (b) edge (n2) (b) edge (n3);
			\path (c) edge (n2) (c) edge (n3);
			\node at (-2.2,-0.75) {$S$};
		\end{scope}
		
		\begin{scope}[xshift=4.5cm, yshift=0.35cm]
			\node[I] (Fa) at (0,-0.971)   {}; \node[below left]  at (Fa){$1$};
			\node[I] (Fb) at (0.9,0.379)  {}; \node[above left] at (Fb){$2$};
			\node[I] (Fc) at (2.1,-0.971) {}; \node[below right] at (Fc){$3$};
			\draw[Fedge] (Fa) -- (Fb) node[midway,left]  {$4$};
			\draw[Fedge] (Fa) -- (Fc) node[midway,below] {$6$};
			\node at (2,0.4) {$\Phi(S)$};
		\end{scope}
		
	\end{tikzpicture}
	\caption{$\deg(S) = \sigma_{12} + \sigma_{13} + \sigma_{23} = 4 + 6 + 0 = 10$}
	\label{ejemplo.cantidad.P4.ind.split}
\end{figure}
The whole of this section rests on the following lemma, which says that two induced $P_4$'s of $S$ sharing a vertex cannot have their endpoints in different components of $\Phi$.

\begin{lemma}
	\label{lema.P4s.que.se.cortan}
	Let $(S,K,I)$ be a split graph and let $H$, $H'$ be induced $P_4$'s of $S$ such that $V(H)\cap V(H')\neq\varnothing$. Then the endpoints of $H$ and the endpoints of $H'$ all lie in the same component of $\Phi(S,K,I)$.
\end{lemma}

\begin{proof}
	Let $u,u'$ be the endpoints of $H$ and let $v,v'$ be the endpoints of $H'$; all four lie in $I$, and $uu',vv'\in E(\Phi)$, where $\Phi=\Phi(S,K,I)$. Since the endpoints of an induced $P_4$ of a split graph lie in $I$ while its internal vertices lie in $K$, any vertex of $V(H)\cap V(H')$ is either an endpoint of both or an internal vertex of both.
	
	In the first case the two edges $uu'$ and $vv'$ of $\Phi$ share a vertex, and the conclusion is immediate. So assume that $V(H)\cap V(H')$ contains an internal vertex $x\in K$ of both. Being internal to $H$, $x$ is adjacent to exactly one endpoint of $H$, so $x\in N_u\bigtriangleup N_{u'}$; likewise $x\in N_v\bigtriangleup N_{v'}$. As the statement is symmetric in $u\leftrightarrow u'$ and in $v\leftrightarrow v'$, we may assume
	\[ x\in (N_u-N_{u'})\cap(N_v-N_{v'}). \]
	Suppose, for contradiction, that $u$ and $v$ lie in different components of $\Phi$. Then $u'$ lies in the component of $u$ and $v'$ in that of $v$, so by \Cref{prop.incomparabilidad}(2) each of the pairs $\{N_u,N_v\}$, $\{N_u,N_{v'}\}$, $\{N_{u'},N_v\}$ and $\{N_{u'},N_{v'}\}$ is comparable. From $x\in N_u-N_{v'}$ we get $N_u\not\subseteq N_{v'}$, hence $N_{v'}\subsetneq N_u$; symmetrically, $x\in N_v-N_{u'}$ gives $N_{u'}\subsetneq N_v$. Now $uu'\in E(\Phi)$ implies that $N_u$ and $N_{u'}$ are incomparable, so there is $y\in N_{u'}-N_u$, and $N_{u'}\subseteq N_v$ gives $y\in N_v-N_u$. Likewise $vv'\in E(\Phi)$ yields some $z\in N_{v'}-N_v$, and $N_{v'}\subseteq N_u$ gives $z\in N_u-N_v$. Therefore $N_u$ and $N_v$ are incomparable, so $uv\in E(\Phi)$ by \Cref{prop.incomparabilidad}(2), contradicting that $u$ and $v$ lie in different components.
\end{proof}

\begin{theorem}
	\label{A_4.conexo.implica.Phi.conexo}
	Let $(S,K,I)$ be a split graph with $\deg(S)\geq 1$. If $A_4(S)$ is connected, then $S$ is active and $\Phi(S)$ is connected.
\end{theorem}

\begin{proof}
	The hypothesis implies that $A_4(S)$ does not have isolated vertices, which means that $S$ is active. Let $u,v\in I$; we may assume $u\neq v$. Since $A_4(S)$ is connected, there is a path $w_0w_1\cdots w_n$ in $A_4(S)$ with $n\geq 1, w_0=u$ and $w_n=v$. For each $i\in[n]$, the edge $w_{i-1}w_i$ of $A_4(S)$ comes from some $H_i\approx P_4$ induced in $S$ with $w_{i-1},w_i\in V(H_i)$. For $1\leq i<n$ we have $w_i\in V(H_i)\cap V(H_{i+1})$, so \Cref{lema.P4s.que.se.cortan} places the endpoints of $H_i$ and those of $H_{i+1}$ in a common component of $\Phi(S)$. By transitivity, the endpoints of $H_1,\ldots,H_n$ all lie in a single component $\Psi$ of $\Phi(S)$. Finally, $u\in V(H_1)\cap I$ and $v\in V(H_n)\cap I$, so $u$ and $v$ are endpoints of $H_1$ and $H_n$ respectively; hence $u,v\in V(\Psi)$.
\end{proof}

\begin{corollary}
	\label{S.indecomp.implies.PhiS.conn}
	If $(S,K,I)$ is a split graph with $|I|\geq 2$ and $\Phi(S,K,I)$ is disconnected, then $S$ is decomposable.
\end{corollary}

\begin{proof}
	If $\deg(S)=0$, both $A_4(S)$ and $\Phi(S,K,I)$ are edgeless. If also $|I|\geq 2$, they are both edgeless and disconnected, and so $S$ is decomposable by \Cref{indecomp.characterization}. If $\Phi(S,K,I)$ is disconnected and $\deg(S)\geq 1$, then \Cref{A_4.conexo.implica.Phi.conexo} implies that $A_4(S)$ is disconnected. Thus, $S$ is decomposable by \Cref{indecomp.characterization}.
\end{proof}

The next result shows that the converse of \Cref{A_4.conexo.implica.Phi.conexo} also holds. As we can see from \Cref{phi.conexo.pero.A_4.no}, the connectedness of $\Phi(S,K,I)$ alone does not necessarily imply the connectedness of $A_4(S)$. In fact, a careful look at the (balanced) split graph $S$ in \Cref{phi.conexo.pero.A_4.no} reveals that it is not active, because the vertex $4$ is universal (see \cite{pastine20252}, page 5). 
\begin{figure}[h]
	\centering
	\begin{tikzpicture}[
		scale=1,
		K/.style={circle,draw,fill=black,minimum size=5.5pt,inner sep=0pt},
		I/.style={circle,draw,fill=white,minimum size=5.5pt,inner sep=0pt},
		every edge/.style={draw,line width=.5pt}
		]
		
		\begin{scope}
			\node[K] (a) at (0,0)       {}; \node[below left] at (a) {$4$};
			\node[K] (b) at (0.65,1.15) {}; \node[above right]      at (b) {$5$};
			\node[K] (c) at (1.6,0.2)   {}; \node[below right]      at (c) {$6$};
			\node[I] (u1) at (-1.55,0.3){}; \node[left]       at (u1){$1$};
			\node[I] (u2) at (-0.9,1.4) {}; \node[below left] at (u2){$2$};
			\node[I] (u3) at (0.6,-1.2) {}; \node[below right]at (u3){$3$};
			\path (a) edge (b) (b) edge (c) (a) edge (c);
			\path (u1) edge (a) (u1) edge (b);
			\path (u2) edge (a) (u2) edge (b);
			\path (u3) edge (a) (u3) edge (c);
			\node at (0.1,-2) {$S$};
		\end{scope}
		
		\begin{scope}[xshift=4.6cm]
			\node[I] (A1) at (-0.75,1.1)  {}; \node[above left]  at (A1){$1$};
			\node[I] (A2) at ( 0.75,1.1)  {}; \node[above right] at (A2){$2$};
			\node[K] (Ab) at ( 1.25,-0.15){}; \node[below right]       at (Ab){$5$};
			\node[I] (A3) at ( 0,-1.15)   {}; \node[below right]       at (A3){$3$};
			\node[K] (Ac) at (-1.25,-0.15){}; \node[below left]        at (Ac){$6$};
			\node[K] (Aa) at ( -1,-1)  {}; \node[below left]       at (Aa){$4$};
			\path (A1) edge (A3) (A1) edge (Ab) (A1) edge (Ac);
			\path (A2) edge (A3) (A2) edge (Ab) (A2) edge (Ac);
			\path (A3) edge (Ab) (A3) edge (Ac);
			\path (Ab) edge (Ac);
			\node at (0,-2) {$A_4(S)$};
		\end{scope}
		
		\begin{scope}[xshift=8.0cm,yshift=-0.7cm]
			\node[I] (P1) at (0,1.15)   {}; \node[above left]  at (P1){$1$};
			\node[I] (P2) at (1.5,1.15) {}; \node[above left] at (P2){$2$};
			\node[I] (P3) at (0,-0.35)  {}; \node[below left]  at (P3){$3$};
			\path (P1) edge (P3);
			\path (P2) edge (P3);
			\node at (0.75,-1.3) {$\Phi(S)$};
		\end{scope}
		
	\end{tikzpicture}
	\caption{$\Phi(S)$ is connected but $A_4(S)$ is not.}
	\label{phi.conexo.pero.A_4.no}
\end{figure} 

\begin{theorem}
	\label{thm.actividad_y_PhiS.conn}
	If $S$ is an active split graph and $\Phi(S)$ is connected, then $S$ is indecomposable.
\end{theorem}

\begin{proof}
	Ignoring parallel edges in $\Phi(S)$, it is clear that $\Phi(S)$ is a subgraph of $A_4(S)$. Furthermore, every $K$-vertex must be connected in $A_4(S)$ to some $I$-vertex, since $S$ is active. Hence, $A_4(S)$ is connected. Applying \Cref{indecomp.characterization}, we conclude that $S$ is indecomposable.
\end{proof}

\begin{theorem}
	\label{S.primo.iff.Phi(S).conexo}
	Let $S$ be an active split graph. Then: 
	\begin{enumerate}
		\item $S$ is prime if and only if $\Phi(S)$ is connected;
		
		\item if $S_k$ is a prime split graph for each $k\in[n]$ and $S=S_n\circ\ldots\circ S_1$, then $\Phi(S)=\dot{\bigcup}_{k=1}^n\Phi(S_k)$;
	\end{enumerate}
\end{theorem}

\begin{proof}
	\begin{enumerate}[(1).]
		\item If $\Phi(S)$ is connected, then $S$ is prime by \Cref{thm.actividad_y_PhiS.conn}. Conversely, if $S$ is prime, then $A_4(S)$ is connected by \Cref{indecomp.characterization}. Hence, $\Phi(S)$ is connected by \Cref{A_4.conexo.implica.Phi.conexo}.
		
		\item It suffices to consider $n=2$, since the general case follows by induction. Let $(K,I)$ be the bipartition for $S$ and let $\Phi=\Phi(S)$. If $(S_1,K_1,I_1)$ and $(S_2,K_2,I_2)$ are prime, then $\Phi(S_1)$ and $\Phi(S_2)$ are connected by (1). If $u,v\in I_1\subseteq I_1\dot{\cup} I_2=I$, then $\sigma_{S_1}(uv)=\sigma_S(uv)$ since $N_{S_1}(w)=N_S(w)$ for all $w\in I_1$. Similarly, $\sigma_{S_2}(uv)=\sigma_S(uv)$ for all $u,v\in I_2$. Hence, $\Phi(S_1)=\Phi[I_1]$ and $\Phi(S_2)=\Phi[I_2]$. If $u\in I_1$ and $v\in I_2$, then $\sigma_S(uv)=0$ by \Cref{prop.incomparabilidad}, because $N_S(u)\subseteq K_1\subseteq N_S(v)$. Therefore, $\Phi(S_1)\dot{\cup}\Phi(S_2)=\Phi[I_1]\dot{\cup}\Phi[I_2]=\Phi[I_1\cup I_2]=\Phi$. \qedhere
	\end{enumerate}
\end{proof}



\section{Multiplicities, neighborhoods and homogeneity} \label{sec:homogen_split}

The first result of this section lists four key properties of $\sigma_{uv}$, relating it to the neighborhoods in $S$ and $\Phi$ of $u$ and $v$. In this context, the concept of twin vertices naturally arises.
Two vertices $u$ and $v$ of a multigraph $M$ are said to be \emph{twins} in $M$ if the following equality of multisets holds: $N_M(u)-m_v\{v\}=N_M(v)-m_u\{u\}$, where $m_u$ is the multiplicity of $u$ in $N_M(v)$, and $m_v$ is the multiplicity of $v$ in $N_M(u)$. It is easy to see that ``to be twins'' is an equivalence relation on $V(M)$. 

\begin{proposition}
	\label{prop.basicas.sigma_uv}
	If $(S,K,I)$ is a split graph, let $\Phi=\Phi(S,K,I)$. Then:
	\begin{enumerate}
		\item $\sigma_{uv}=0$ (i.e., $uv\notin E(\Phi)$) and $d_v \leq d_u$ if and only if $N_v \subseteq N_u$; 
		
		\item $\sigma_{uv}=0$ and $d_u =d_v$ if and only if $N_u = N_v$;
		
		\item if $N_u = N_v$, then $\sigma_{uw}=\sigma_{vw}$ for every $w\in I-\{u,v\}$; in other words, if $u$ and $v$ are twins in $S$, then they are also twins in $\Phi$;
		
		\item if $\sigma_{uv}=1$, then $d_u =d_v$ and $|N_u -N_v |=|N_v -N_u |=1$;
	\end{enumerate}
\end{proposition}

\begin{proof}
	\begin{enumerate}[(1).]
		\item If $\sigma_{uv}=|N_u-N_v||N_v-N_u|=0$, then $|N_u-N_v|=0$ or $|N_v-N_u|=0$. Thus, $N_u\subseteq N_v$ or $N_v\subseteq N_u$.
		
		\item It follows from (1).
		
		\item Assume $N_u=N_v$ and let $w\in I-\{u,v\}$. Then $N_u-N_w=N_v-N_w$ and $N_w-N_u=N_w-N_v$, so $\sigma_{uw}=|N_u-N_w|\,|N_w-N_u|=|N_v-N_w|\,|N_w-N_v|=\sigma_{vw}$. Hence $\sigma_{uw}=\sigma_{vw}$ for every $w\in I-\{u,v\}$, so $u$ and $v$ have identical neighborhoods in $\Phi$ counting multiplicities; that is, they are twins in $\Phi$.
		
		\item Since $1=\sigma_{uv}=|N_u-N_v||N_v-N_u|$, we have $|N_u-N_v|=|N_v-N_u|=1$. Also, $d_u-\eta_{uv}=|N_u-N_v|=|N_v-N_u|=d_v-\eta_{uv}$ implies $d_u=d_v$. \qedhere
	\end{enumerate}
\end{proof}

Notice that the converse of \Cref{prop.basicas.sigma_uv}(3) does not hold in general. As a counterexample, let $(S,K,I)$ be the (active) split graph with $I=[4]$, $K=\{5,6,7,8\}$, $N_1=N_2=\{5\}$, $N_3=\{6,7\}$ and $N_4=\{6,8\}$. Then, $\sigma_{13}=\sigma_{14}=\sigma_{23}=\sigma_{24}=2$, $\sigma_{34}=1$ and $\sigma_{12}=0$. So, $3$ and $4$ are twins in $\Phi$ but not in $S$.

Let $(S,K,I)$ be a split graph. We say that $(S,K,I)$ is $\delta$-\emph{homogeneous} if $\deg_S(v)=\delta$ for all $v\in I$. 

\begin{proposition}
	\label{S.homogeneo.implica...}
	Let $(S,K,I)$ be a $\delta$-homogeneous split graph with $\deg(S)\geq 1$, and let $\Phi=\Phi(S,K,I)$. Then:
	\begin{enumerate}
		\item $\sigma_{uv}=0$ if and only if $u$ and $v$ are twins in $S$;
		
		\item $\delta\notin\{0,|K|\}$;
		
		\item $\Phi$ cannot contain 3 vertices $u,v,w$ with $\sigma_{uv}\neq 0$ and $\sigma_{uw}=0=\sigma_{vw}$;
		
		\item $\text{diam}(\Phi)\leq 2$;
	\end{enumerate}
\end{proposition}

\begin{proof}
	\begin{enumerate}[(1).]
		\item It follows immediately from \Cref{prop.basicas.sigma_uv}(2).
		
		\item If $\delta \in\{0,|K|\}$, then $N_v\in\{\varnothing, K\}$ for all $v\in I$. So, $\Phi$ is clearly edgeless, contradicting the hypothesis.
	    
	    \item Suppose we have $\{u,v,w\}\subseteq I$ with $\sigma_{uv}\neq 0$ and $\sigma_{uw}=0=\sigma_{vw}$. By \Cref{prop.basicas.sigma_uv}(2) we have $N_u =N_w$ and $N_w =N_v$. Thus $N_u =N_v$, and hence $\sigma_{uv}=0$ by \Cref{prop.basicas.sigma_uv}(2), a contradiction.
		
		\item If $\text{diam}(\Phi)\ge 3$, it is easy to see that we contradict (3).
	\end{enumerate}
\end{proof}

\begin{theorem}
	\label{Phi.simple.conexo.1}
	Let $(S,K,I)$ be a split graph such that $K=\bigcup_{v\in I}N_v$ and $|I|\geq 2$. If $\Phi(S,K,I)$ is simple and connected, then $(S,K,I)$ is $\delta$-homogeneous and balanced.
\end{theorem}

\begin{proof}
	Since $\Phi=\Phi(S,K,I)$ is simple, $uv\in E(\Phi)$ if and only if $\sigma_{uv}=1$. Thus, whenever $uv,vw\in E(\Phi)$, \Cref{prop.basicas.sigma_uv}(4) gives $d_u=d_v$ and $d_v=d_w$. As $\Phi$ is connected, all vertices of $I$ share a common degree $\delta$ in $S$; that is, $S$ is $\delta$-homogeneous. Since the hypothesis clearly imply $\deg(S)\ge 1$, we have by \Cref{S.homogeneo.implica...}(2) that $0<\delta<|K|$. 
	
	It remains to prove that $S$ is balanced, i.e. $\omega(S)=|K|$ and $\alpha(S)=|I|$. Since $I$ is an independent set, every clique of $S$ contains at most one vertex of $I$, so $\omega(S)\leq|K|+1$, with $\omega(S)=|K|+1$ only if some $u\in I$ is adjacent to every vertex of $K$, i.e. $d_u=|K|$. But $d_u=\delta<|K|$ for all $u\in I$, so $\omega(S)=|K|$. Since $K$ is a clique, every independent set of $S$ contains at most one vertex of $K$, so $\alpha(S)\leq|I|+1$, with $\alpha(S)=|I|+1$ only if some $y\in K$ has no neighbor in $I$. But $K=\bigcup_{v\in I}N_v$ means every $x\in K$ has a neighbor in $I$, so $\alpha(S)=|I|$.
\end{proof}

\section{Simple and complete factor graphs}
\label{sec:simple_complete_Phi}


%
To prove \Cref{Phi.simple.completo.implica...}, we need  \Cref{intersec_family}, an auxiliary result about certain intersecting families of sets. \Cref{intersec_family} is a translation to our language and notation of Lemma 1 in \cite{automorphism.johnson.graph}, and is closely related to the description of the cliques of a Johnson graph given in \cite{shuldiner.oldford.johnson.cliques}.

\begin{theorem}[\cite{automorphism.johnson.graph}, page 267]
	\label{intersec_family}
	Let $|I|\geq 3$ and let $\{N_v :v\in I\}$ be a family of finite sets of size $\delta\geq 1$ such that $|N_u \cap N_v |=\delta-1$ for all $\{u,v\}\subseteq I$. Let also $K=\bigcup_{v\in I}N_v$. 
	\begin{enumerate}[1.]
		\item If $B$ is a non-empty subset of $I$, then 
		\begin{equation*}
			\left|\bigcap_{v\in B} N_v \right| \in\{\delta+1-|B|,\delta-1\}.
		\end{equation*}
		
		\item The following 3 statements are equivalent:
		\begin{enumerate}[2.1.]
			\item there exists a triple $T\subseteq I$ such that $|\bigcap_{v\in T}N_v|=\delta-1$;
			\item \[ \left|\bigcap_{v\in A}N_v \right|=\delta-1, \] for all $A\subseteq I$ with $|A|\geq 2$;
			\item \[ |K|=|I|+\delta-1. \] 
		\end{enumerate}
		
		\item The following 3 statements are equivalent:
		\begin{enumerate}[3.1.]
			\item there exists a triple $T\subseteq I$ such that $|\bigcap_{v\in T}N_v|=\delta-2$;
			\item \[ \left|\bigcap_{v\in A}N_v \right|=\delta+1-|A|, \] for all $A\subseteq I$ with $A\neq\varnothing$;
			\item \[ |K|=\delta+1. \] 
		\end{enumerate}
	\end{enumerate}
\end{theorem}

\begin{corollary}
	\label{intersect_family_corol}
	Let $\{N_v :v\in I, |I|\geq 2\}$ be a family of finite sets of size $\delta\geq 1$ such that $|N_u \cap N_v |=\delta-1$ for all $\{u,v\}\subseteq I$. Let also $K=\bigcup_{v\in I}N_v$ and $U=\bigcap_{v\in I}N_v$. Then:
	\begin{enumerate}
		\item $|U|\in\{\delta-1,\ \delta+1-|I|\}$;
		
		\item $|U|=\delta-1$ if and only if $|K|=|I| +\delta-1$;
		
		\item $|U|=\delta+1-|I|$ if and only if $|K| =\delta+1$;
		
		\item $|K| =|I| +|U|$ (in particular, $|K|\geq |I|$).
	\end{enumerate}
\end{corollary}

\begin{proof}
	Let $\omega=|K|$ and $\alpha=|I|$.
	\begin{enumerate}[(1).]
		\item If $\alpha=2$, the claim is immediate. If $\alpha\geq 3$, apply \Cref{intersec_family}(1) to $B=I$.
		
		\item If $\alpha=2$, we know that $|U|=\delta-1$; moreover $\omega=|N_u\cup N_v|=2\delta-(\delta-1)=\delta+1=\alpha+\delta-1$. Hence both sides of the equivalence always hold, and there is nothing to prove.
		
		Assume now $\alpha\geq 3$. If $\omega=\alpha+\delta-1$, then \Cref{intersec_family}(2) gives $|\bigcap_{v\in A}N_v|=\delta-1$ for every $A\subseteq I$ with $|A|\geq 2$; taking $A=I$ yields $|U|=\delta-1$. Conversely, assume $|U|=\delta-1$ and let $T=\{t_1,t_2,t_3\}\subseteq I$ be any triple. Then $U\subseteq\bigcap_{v\in T}N_v\subseteq N_{t_1}\cap N_{t_2}$, and $|U|=\delta-1=|N_{t_1}\cap N_{t_2}|$, so $|\bigcap_{v\in T}N_v|=\delta-1$. Now \Cref{intersec_family}(2) gives $\omega=\alpha+\delta-1$.
		
		\item If $\alpha=2$, both sides always hold, by the computation in (2) and because $\delta+1-\alpha=\delta-1$.
		
		Assume $\alpha\geq 3$. If $\omega=\delta+1$, then \Cref{intersec_family}(3) gives $|\bigcap_{v\in A}N_v|=\delta+1-|A|$ for every non-empty $A\subseteq I$; taking $A=I$ yields $|U|=\delta+1-\alpha$. Conversely, assume $|U|=\delta+1-\alpha$. Since $\alpha\geq 3$, we have $|U|<\delta-1$, so $\omega\neq \alpha+\delta-1$, by (2). Thus, no triple $T\subseteq I$ satisfies $|\bigcap_{v\in T}N_v|=\delta-1$, by \Cref{intersec_family}(2). But \Cref{intersec_family}(1) applied to any triple $T$ gives $|\bigcap_{v\in T}N_v|\in\{\delta-2,\delta-1\}$, whence $|\bigcap_{v\in T}N_v|=\delta-2$. Now \Cref{intersec_family}(3) gives $\omega=\delta+1$.
		
		\item By (1), either $|U|=\delta-1$ or $|U|=\delta+1-\alpha$. In the first case, $\omega=\alpha +\delta-1=\alpha+|U|$, by (2). In the second case, $\omega=\delta+1=(\delta+1-\alpha)+\alpha=|U|+\alpha$, by (3). \qedhere
	\end{enumerate}
\end{proof}

When $(S,K,I)$ is a split graph with $|I|\geq 2$ and $\Phi(S,K,I)$ is simple and complete, \Cref{prop.basicas.sigma_uv}(4) gives $|N_S(v)|=\delta$ for all $v\in I$ and $|N_S(u)\cap N_S(v)|=\delta-1$ for all $\{u,v\}\subseteq I$, for a common $\delta\geq 1$. Thus, \Cref{intersect_family_corol} applies to $\{N_S(v):v\in I\}$, and this is the source of items (3),...,(6) of \Cref{Phi.simple.completo.implica...} below. Before proving this theorem, we need a technical lemma.

\begin{lemma}
	\label{Phi.completo.implica.inactivos.universales}
	Let $(S,K,I)$ be a balanced split graph such that $|I|\geq 2$ and $\Phi(S)$ is complete. If $x$ is an inactive vertex in $S$, then $x$ is universal.
\end{lemma}

\begin{proof}
	The assumptions on $S$ guarantee $I\subseteq \act(S)$. Hence, if $x\notin \act(S)$, then $x\in K$. Suppose $x$ is not universal. Since $S$ has no swing vertices, it follows that $N_x\cap I\neq\varnothing$. Then there exist $u,v\in I$ such that $u\in N_x$ and $v\notin N_x$. On the other hand, since $\sigma_{uv}>0$, there exist $y,z\in K$ such that $uyzv\prec S$. But then $uxzv\prec S$, contradicting that $x$ is inactive.
\end{proof}

\begin{theorem}
	\label{Phi.simple.completo.implica...}
	Let $(S,K,I)$ be a split graph such that $\bigcup_{v\in I}N_v = K$, $|I|\geq 2$ and $\Phi(S,K,I)$ is simple and complete. Let $U$ be the set of universal vertices of $S$, $\omega=\omega(S)$ and $\alpha=\alpha(S)$. Then:
	\begin{enumerate}
		\item $(S,K,I)$ is $\delta$-homogeneous and balanced; in particular, $(K,I)$ is the unique bipartition for $S$ and $(|K|,|I|)=(\omega,\alpha)$;
		
		\item $1\leq \delta\leq \omega-1$;
		
		\item $U=\bigcap_{v\in I}N_v$ and $|U|\in\{\delta-1,\ \delta+1-\alpha\}$;
		
		\item $|U|=\delta-1$ if and only if $\omega=\alpha +\delta-1$;
		
		\item $|U|=\delta+1-\alpha$ if and only if $\omega =\delta+1$;
		
		\item $\omega =\alpha +|U|$ (in particular, $\omega\geq \alpha$);    
		
		\item $S$ is active if and only if $U=\varnothing$;
		
		\item if $S$ is active, then $\omega=\alpha$ and $\delta\in\{1,\omega-1\}$;
		
		\item if $\omega=\alpha$ or $\delta=1$, then $S$ is active.
	\end{enumerate} 
\end{theorem}

\begin{proof}
	\begin{enumerate}[(1).]
		\item It follows immediately from \Cref{Phi.simple.conexo.1}.
		
		\item It follows immediately from (1) and \Cref{S.homogeneo.implica...}(2).
		
		\item We know that $U=\bigcap_{v\in I}N_v$ from Proposition 4.1(1) of \cite{pastine20252}. The rest of the statement follows from \Cref{intersect_family_corol}(1).
		
		\item It is the content of \Cref{intersect_family_corol}(2).
		
		\item It is the content of \Cref{intersect_family_corol}(3).
		
		\item It is the content of \Cref{intersect_family_corol}(4).
		
		\item If $S$ is active, then $U=\varnothing$, since universal vertices are inactive (see \cite{pastine20252}, page 5). Conversely, assume $U=\varnothing$. By (1), $S$ is balanced, so \Cref{Phi.completo.implica.inactivos.universales} applies and every inactive vertex of $S$ is universal. As $U=\varnothing$, the graph $S$ has no inactive vertices; that is, $S$ is active.
		
		\item If $S$ is active, then $U=\varnothing$ by (7), and hence $\omega=\alpha$ by (6). Moreover, $0=|U|\in\{\delta-1,\delta+1-\alpha\}$ by (3), so $\delta=1$ or $\delta=\alpha-1=\omega-1$.
		
		\item If $|K|=|I|$, then $|U|=0$ by (6), and $S$ is active by (7). If $\delta=1$, then $|U|\in\{0,2-|I|\}$ by (3). Since $|U|\geq 0$ and $|I|\geq 2$, the second alternative forces $|I|=2$ and $|U|=0$. In either case $U=\varnothing$, so $S$ is active by (7). \qedhere
	\end{enumerate} 
\end{proof}

Regarding \Cref{Phi.simple.completo.implica...}(9), it is important to note that $\delta=\omega-1$ does not generally imply that $S$ is active. For example, consider the split graph $(S,K,I)$ with $I=[2]$, $K=\{3,4,5,6\}$, where $N_1=\{3,4,5\}$ and $N_2=\{4,5,6\}$. Here, we see that $d_1=3=\omega-1=d_2$ and $\sigma_{12}=1$, so $\Phi(S)$ is simple and complete. However, $S$ is not active, since vertices $4$ and $5$ are universal.

\begin{theorem}
	\label{act.Phi.simple.compl.iff...}
	Let $(S,K,I)$ be an active split graph. Then $\Phi(S)$ is simple and complete if and only if $|K|=|I|$ and one of the following holds:
	\begin{enumerate}
		\item the edges between $K$ and $I$ form a perfect matching; that is, there is a bijection $v\mapsto x_v$ from $I$ to $K$ with $N_S(v)=\{x_v\}$ for all $v\in I$ (every vertex of $I$ has a unique neighbor in $K$); or
		
		\item the edges between $K$ and $I$ form the complement of a perfect matching; that is, there is a bijection $v\mapsto x_v$ from $I$ to $K$ with $N_S(v)=K-\{x_v\}$ for all $v\in I$ (every vertex of $I$ has a unique non-neighbor in $K$).
	\end{enumerate}
\end{theorem}

\begin{proof}
	($\Leftarrow$). In case (1), $N_u$ and $N_v$ are distinct singletons for $u\neq v$, so $\sigma_{uv}=1$ for every $\{u,v\}\subseteq I$, by \Cref{prop.incomparabilidad}(1). Hence, $\Phi=\Phi(S)$ is simple and complete. Each $v\in I$ is an endpoint of the induced path $v\,x_v\,x_u\,u$ for any $u\neq v$, so $S$ is active. In case (2), for $u\neq v$ we have $N_u-N_v=\{x_v\}$ and $N_v-N_u=\{x_u\}$, so again $\sigma_{uv}=1$ by \Cref{prop.incomparabilidad}(1), and $\Phi$ is simple and complete. Moreover, $vx_ux_vu\prec S$, showing that $S$ is active. 
	
	($\Rightarrow$). By \Cref{Phi.simple.completo.implica...}(8), $|K|=|I|$ and $\delta\in\{1,\,|K|-1\}$. If $\delta=1$, each $N_v$ is a singleton; since $\Phi$ is complete, $N_u\neq N_v$ for $u\neq v$, so the singletons are distinct and $v\mapsto x_v$ (where $N_v=\{x_v\}$) is a bijection $I\to K$. This is case (1). If $\delta=|K|-1$, each $N_v$ is $K$ minus a single vertex $x_v$. Suppose that $u\neq v$ but $x_u=x_v$. Then $N_u=K-\{x_u\}=K-\{x_v\}=N_v$, so $\sigma_{uv}=0$; but this contradicts the completeness of $\Phi$. Hence $u\mapsto x_u$ is injective, and since $|K|=|I|$, it is a bijection $I\to K$ with $N_v=K-\{x_v\}$ for every $v\in I$. This is case (2).
\end{proof}


\section{Simple and connected factor graphs}
\label{sec:simple_Phi}


The following is a technical result that relates the subgraph of $\Phi$, induced by a subset $A$ of its vertices, to the factor graph of the corresponding induced subgraph in $S$ that has $A$ as an independent part. It also shows that, in some way, $\Phi$ does not depend on universal vertices, i.e., removing universal vertices in the underlying split graph does not alter the structure of the factor graph.

\begin{proposition}
	\label{subgrafos.de.Phi}
	Let $(S,K,I)$ be a split graph and let $A\subseteq I$ with $|A|\geq 2$. Write $K_A=\bigcup_{v\in A}N_v$ and $K_A'=K_A-\bigcap_{v\in A}N_v$, and define
	\[ S_A=S[A\cup K_A], \qquad S_A'=S_A-\bigcap_{v\in A}N_v. \]
	Then $(K_A,A)$ and $(K_A',A)$ are bipartitions of $S_A$ and $S_A'$, respectively, and
	\[ \Phi(S_A,K_A,A)=\Phi(S,K,I)[A]=\Phi(S_A',K_A',A). \]
\end{proposition}
\begin{proof}
	By construction, it is clear that $(K_A,A)$ and $(K_A',A)$ are bipartitions of $S_A$ and $S_A'$, respectively. All three multigraphs have vertex set $A$, so it suffices to show that $\sigma_{S_A}(uv)=\sigma_{uv}=\sigma_{S_A'}(uv)$, for every pair $u,v\in A$. Recall from \Cref{prop.incomparabilidad}(1) that $\sigma_{uv}=(d_u-\eta_{uv})(d_v-\eta_{uv})$. Since $N_w\subseteq\bigcup_{w\in A}N_w$ for every $w\in A$, we have $N_{S_A}(w)=N_w$, so $\deg_{S_A}(w)=d_w$ and $\eta_{uv}(S_A)=\eta_{uv}$; hence $\sigma_{S_A}(uv)=\sigma_{uv}$. For $S_A'$, write $\eta_A=\big|\bigcap_{w\in A}N_w\big|$, $d_w'=\deg_{S_A'}(w)$ and $\eta_{uv}'=\eta_{uv}(S_A')$. By construction,
	\[ N_{S_A'}(w)=N_w-\bigcap_{a\in A}N_a \quad (w\in A), \quad N_{S_A'}(u)\cap N_{S_A'}(v)=(N_u\cap N_v)-\bigcap_{a\in A}N_a. \]
	Since $w\in A$, we have $\bigcap_{a\in A}N_a\subseteq N_w$, and likewise $\bigcap_{a\in A}N_a\subseteq N_u\cap N_v$; therefore
	\[ d_w'=d_w-\eta_A \quad (w\in A), \qquad \eta_{uv}'=\eta_{uv}-\eta_A. \]
	Consequently $d_u'-\eta_{uv}'=(d_u-\eta_A)-(\eta_{uv}-\eta_A)=d_u-\eta_{uv}$, and similarly $d_v'-\eta_{uv}'=d_v-\eta_{uv}$. Hence
	\[ \sigma_{S_A'}(uv)=(d_u'-\eta_{uv}')(d_v'-\eta_{uv}')=(d_u-\eta_{uv})(d_v-\eta_{uv})=\sigma_{uv}. \qedhere \]
\end{proof} 

The following lemma is needed to prove \Cref{caract.Phi.simple.conexo}, which characterizes the active split graphs having a simple and connected factor graph.

\begin{lemma}
	\label{twin.reduction.lemma}
	Let $(S,K,I)$ be an active split graph such that $\Phi=\Phi(S)$ is simple and connected. Let $J\subseteq I$ contain exactly one vertex from each class of the equivalence relation on $I$ given by $N_u=N_v$, and let $H=S[J\cup K]$. Then:
	\begin{enumerate}
		\item $\deg(S)\geq 1$, and $(S,K,I)$ is balanced and $\delta$-homogeneous with $0<\delta<|K|$; in particular, $\omega(S)=|K|\geq 2$ and $\alpha(S)=|I|\geq 2$;
		
		\item $\bigcup_{v\in I}N_v=K$ and $\bigcap_{v\in I}N_v=\varnothing$;
		
		\item for distinct $u,v\in I$, $\sigma_{uv}=0$ if and only if $N_u=N_v$; hence the classes above are exactly the classes of twins of $S$ that are contained in $I$, and $uv\in E(\Phi)$ if and only if $u$ and $v$ lie in different classes;
		
		\item $N_H(v)=N_S(v)$ for every $v\in J$;
		
		\item $\bigcup_{v\in J}N_v=K$ and $\bigcap_{v\in J}N_v=\varnothing$;
		
		\item $\Phi(H,K,J)=\Phi[J]$, and this multigraph is simple and complete;
		
		\item $H$ is active and $|K|=|J|$.
	\end{enumerate}
\end{lemma}

\begin{proof}
	\begin{enumerate}[(1).]
		\item It easily follows by combining Proposition 4.1(7) of \cite{pastine20252}, Proposition 4.2 of \cite{pastine20252}, \Cref{Phi.simple.conexo.1} and \Cref{S.homogeneo.implica...}(2).
		
		\item The first equality was implicitly obtained in (1). For the second, universal vertices are inactive (see \cite{pastine20252}, page 5), so the set $U$ of universal vertices of $S$ is empty; and $U=\bigcap_{v\in I}N_v$, by Proposition 4.1(1) of \cite{pastine20252}.
		
		\item By (1) we may apply \Cref{S.homogeneo.implica...}(1), which states that $\sigma_{uv}=0$ if and only if $u$ and $v$ are twins in $S$; for $u,v\in I$ this means precisely $N_u=N_v$. The last assertion follows because $uv\in E(\Phi)$ if and only if $\sigma_{uv}\neq 0$.
		
		\item 
		$N_H(v)=N_S(v)\cap(J\cup K)=N_S(v)$.
		
		\item It immediately follows from (3) and (2).
		
		\item By (1), $E(\Phi)\neq\varnothing$, so (3) provides two distinct classes and $|J|\geq 2$. By (5), $\bigcup_{v\in J}N_v=K$, so \Cref{subgrafos.de.Phi} applied to $A=J$ yields $K_A=K$, $S_A=H$ and $\Phi(H,K,J)=\Phi[J]$. This multigraph is simple because $\Phi$ is, and it is complete because distinct members of $J$ lie in different classes and hence satisfy $\sigma_{uv}\neq 0$, by (3).
		
		\item By (4) and (5), $\bigcup_{v\in J}N_H(v)=K$; together with $|J|\geq 2$ and (6), this shows that \Cref{Phi.simple.completo.implica...} applies to $(H,K,J)$. The set of universal vertices of $H$ is $\bigcap_{v\in J}N_v=\varnothing$, by (5), so $H$ is active by \Cref{Phi.simple.completo.implica...}(7). Items (8) and (1) of \Cref{Phi.simple.completo.implica...} then give $|K|=|J|$. \qedhere
	\end{enumerate}
\end{proof}


\begin{theorem}
	\label{caract.Phi.simple.conexo}
	Let $(S,K,I)$ be a split graph with $|K|\geq 2$. Then $S$ is active and $\Phi(S)$ is simple and connected if and only if one of the following holds:
	\begin{enumerate}
		\item every $v\in I$ has a unique neighbor $x_v$ in $K$, and every vertex of $K$ has at least one neighbor in $I$; or
		
		\item every $v\in I$ has a unique non-neighbor $x_v$ in $K$, and every vertex of $K$ has at least one non-neighbor in $I$.
	\end{enumerate}
	In either case $|K|\leq|I|$ and $\Phi(S)$ is the complete multipartite graph whose parts are the classes of twins of $I$.
\end{theorem}

\begin{proof}
	($\Leftarrow$) Assume (1) and let $u,v\in I$. Since $N_u=\{x_u\}$ and $N_v=\{x_v\}$, we get $\sigma_{uv}=1$ if $x_u\neq x_v$, and $\sigma_{uv}=0$ otherwise, by \Cref{prop.incomparabilidad}(1). Assume instead (2). Then $N_u-N_v=\{x_v\}$ and $N_v-N_u=\{x_u\}$ whenever $x_u\neq x_v$, and $N_u=N_v$ otherwise; so again $\sigma_{uv}=1$ if $x_u\neq x_v$, and $\sigma_{uv}=0$ otherwise.
	
	In both cases $\Phi(S)$ is simple, and it is precisely the complete multipartite graph with parts $\{v\in I:x_v=x\}$, $x\in K$. Since every vertex of $K$ is of the form $x_v$ for some $v\in I$, there are exactly $|K|\geq 2$ non-empty parts; hence $\Phi(S)$ is connected and $|K|\leq|I|$.
	
	It remains to prove that $S$ is active. Let $v\in I$ and pick $u\in I$ with $x_u\neq x_v$, which exists because $|K|\geq 2$ and every vertex of $K$ is some $x_w$. In case (1), $vx_vx_uu\prec S$; in case (2), $vx_ux_vu\prec S$. As every vertex of $K$ equals $x_w$ for some $w\in I$, each vertex of $S$ occurs in a path of this form, and therefore $\act(S)=V(S)$.
	
	($\Rightarrow$) Let $J$ and $H$ be as in \Cref{twin.reduction.lemma}. By \Cref{act.Phi.simple.compl.iff...} applied to $(H,K,J)$, there is a bijection $v\mapsto x_v$ from $J$ to $K$ such that either $N_H(v)=\{x_v\}$ for all $v\in J$, or $N_H(v)=K-\{x_v\}$ for all $v\in J$. By \Cref{twin.reduction.lemma}(4), $N_S(v)=N_H(v)$ for every $v\in J$; and every $u\in I-J$ is a twin of some $v\in J$ by \Cref{twin.reduction.lemma}(3), so $N_S(u)=N_S(v)$. Defining $x_u=x_v$ for every such $u$ extends the assignment $v\mapsto x_v$ from $J$ to all of $I$; the vertices of $I$ that share the same image $x\in K$ are precisely the members of one twin class. Since $v\mapsto x_v$ was already a bijection from $J$ onto $K$, every $x\in K$ appears as $x_v$ for some $v$, which gives the surjectivity condition: every vertex of $K$ has a neighbor in $I$ in case~(1), respectively a non-neighbor in $I$ in case~(2). \qedhere
\end{proof}

\Cref{act.Phi.simple.compl.iff...} is exactly the twin-free case of \Cref{caract.Phi.simple.conexo}. Indeed, by \Cref{twin.reduction.lemma}(3) two vertices of $I$ are non-adjacent in $\Phi$ precisely when they are twins in $S$, so $\Phi$ is complete if and only if no two vertices of $I$ are twins, if and only if all the parts are singletons, if and only if $|K|=|I|$; and in that case conditions (1) and (2) of \Cref{caract.Phi.simple.conexo} say that the edges between $K$ and $I$ form a perfect matching or the complement of one, which is the content of \Cref{act.Phi.simple.compl.iff...}. In the opposite direction, every graph covered by \Cref{caract.Phi.simple.conexo} is obtained from one covered by \Cref{act.Phi.simple.compl.iff...} by replacing each vertex of $I$ by a class of twins; correspondingly, $\Phi$ passes from $K_{|K|}$ to a complete multipartite blow-up of it. Finally, notice that alternatives (1) and (2) of \Cref{caract.Phi.simple.conexo} exclude each other unless $|K|=2$, in which case they coincide.

\begin{corollary}
	\label{Phi.simple.conexo.deg.cotas}
	Let $(S,K,I)$ be an active split graph with $V(S)\neq\varnothing$, such that $\Phi=\Phi(S)$ is simple and connected. Let $V_1,\dots,V_k$ be the classes of twins of $I$, let $n_i=|V_i|$, and write $\omega=\omega(S)$ and $\alpha=\alpha(S)$. Then:
	\begin{enumerate}
		\item $\omega=\omega(\Phi)=|K|=k\leq|I|=\alpha$, and $\omega=\alpha$ if and only if $\Phi$ is complete, if and only if $n_i=1$ for every $i\in[k]$; in particular, $S$ and $\Phi$ have the same clique number;
		
		\item every vertex of $\Phi$ lies in a maximum clique of $\Phi$, and $\deg_{\Phi}(v)=\alpha-n_i$ for every $v\in V_i$; in particular, $\omega\leq 1+\min\{\deg_{\Phi}(v):v\in I\}$;
		
		\item $\displaystyle\deg(S)=\sum_{i<j}n_in_j=\frac{\alpha^2-\sum_{i=1}^{k}n_i^2}{2}$;
		
		\item $\deg(S)=\binom{\alpha}{2}$ if and only if $\Phi$ is complete.
	\end{enumerate}
\end{corollary}

\begin{proof}
	By \Cref{twin.reduction.lemma}(1), $(S,K,I)$ is balanced with $|K|\geq 2$, so $\omega=|K|$ and $\alpha=|I|$; and by \Cref{twin.reduction.lemma}(3), the classes of twins of $I$ are well defined and two vertices of $I$ are adjacent in $\Phi$ if and only if they lie in different classes. Hence \Cref{caract.Phi.simple.conexo} applies to $(S,K,I)$ and shows that $\Phi$ is the complete multipartite graph with parts $V_1,\dots,V_k$; since these parts are indexed by the vertices of $K$, we have $k=|K|$.
	\begin{enumerate}[(1).]
		\item A clique of a complete multipartite graph meets each part at most once, and picking one vertex from each part produces a clique; therefore $\omega(\Phi)=k=|K|=\omega$. Moreover $\alpha=\sum_{i=1}^{k}n_i\geq k$, with equality if and only if $n_i=1$ for every $i\in[k]$, that is, if and only if all the parts are singletons, which is precisely the condition for $\Phi$ to be complete.
		
		\item Given $v\in V_i$, adding to $v$ one vertex of $V_j$ for each $j\neq i$ yields a clique of size $k$, which is maximum by (1). Since $\Phi$ is the complete multipartite graph with parts $V_1,\dots,V_k$, the vertex $v$ is adjacent exactly to the vertices of $I-V_i$, so $\deg_{\Phi}(v)=\alpha-n_i$. Finally, $\alpha\geq n_i+(k-1)$, whence $\deg_{\Phi}(v)\geq k-1=\omega-1$.
		
		\item The edges of $\Phi$ are exactly the pairs of vertices lying in different parts, so $\deg(S)=|E(\Phi)|=\sum_{i<j}n_in_j$. The second expression follows from $\big(\sum_{i}n_i\big)^2=\sum_{i}n_i^2+2\sum_{i<j}n_in_j$ together with $\sum_{i}n_i=\alpha$.
		
		\item By (3), $\deg(S)=\binom{\alpha}{2}$ if and only if $\alpha^2-\sum_{i}n_i^2=\alpha^2-\alpha$, that is, if and only if $\sum_{i}n_i(n_i-1)=0$. As $n_i\geq 1$ for every $i$, this happens exactly when $n_i=1$ for every $i\in[k]$, which by (1) means that $\Phi$ is complete. \qedhere
	\end{enumerate}
\end{proof}

Finally, we record what happens under complementation. Recall that a graph is active if and only if its complement is. Moreover, if $(S,K,I)$ is a split graph, keep in mind that $\overline{S}$ is split as well, and $(I,K)$ is a bipartition for it. So, $\Phi(\overline{S},I,K)$ is a multigraph with vertex set $K$.  

\begin{corollary}
	\label{Phi.complement.not.simple}
	Let $(S,K,I)$ be an active split graph such that $\Phi(S)$ is simple and connected. Then $\Phi(\overline{S})$ is complete, and it is simple if and only if $\Phi(S)$ is complete. More precisely:
	\begin{enumerate}
		\item if $\Phi(S)$ is complete, then $\Phi(\overline{S})$ is simple and complete;
		
		\item if $\Phi(S)$ is not complete, then $\Phi(\overline{S})$ is complete but not simple.
	\end{enumerate}
\end{corollary}

\begin{proof}
	Write $\omega=|K|$, and let $V_1,\dots,V_k$ and $n_i$ be as in \Cref{Phi.simple.conexo.deg.cotas}, so that $k=\omega$.
	\begin{enumerate}[(1).]
		\item By \Cref{act.Phi.simple.compl.iff...}, the edges of $S$ between $K$ and $I$ form a perfect matching or the complement of one. The edges of $\overline{S}$ between $I$ and $K$ are exactly the non-edges of $S$ between $K$ and $I$, so the complement of a perfect matching becomes a perfect matching and vice versa. Hence \Cref{act.Phi.simple.compl.iff...} applies to $(\overline{S},I,K)$ and shows that $\Phi(\overline{S})$ is simple and complete.
		
		\item Let $J$ and $H$ be as in \Cref{twin.reduction.lemma}. By \Cref{act.Phi.simple.compl.iff...} applied to $(H,K,J)$, the edges of $H$ between $K$ and $J$ form a perfect matching or the complement of one; arguing as in (1), \Cref{act.Phi.simple.compl.iff...} applies to $(\overline{H},J,K)$ and yields that $\Phi(\overline{H},J,K)$ is simple and complete. In particular, $\sigma_{\overline{H}}(xy)\neq 0$ for every $\{x,y\}\subseteq K$. Since $H\prec S$ gives $\overline{H}\prec\overline{S}$, every induced $P_4$ of $\overline{H}$ is an induced $P_4$ of $\overline{S}$, whence $\sigma_{\overline{S}}(xy)\geq\sigma_{\overline{H}}(xy)>0$ for all $\{x,y\}\subseteq K$; that is, $\Phi(\overline{S})$ is complete.
		
		It remains to show that $\Phi(\overline{S})$ is not simple. Since $\Phi(S)$ is not complete, \Cref{Phi.simple.conexo.deg.cotas}(1) provides an index $i$ with $n_i\geq 2$; as $k=\omega\geq 2$, we may pick $j\neq i$, and then \Cref{Phi.simple.conexo.deg.cotas}(3) gives
		\[ \deg(S)=\sum_{r<s}n_rn_s\geq n_in_j+\binom{\omega}{2}-1>\binom{\omega}{2}, \]
		because the sum has $\binom{k}{2}=\binom{\omega}{2}$ terms, each of them at least $1$, and $n_in_j\geq 2$. Since $\deg(\overline{S})=\deg(S)$ by Corollary 6.2 of \cite{pastine20252}, while a simple graph on the $\omega$ vertices of $\Phi(\overline{S})$ has at most $\binom{\omega}{2}$ edges, $\Phi(\overline{S})$ cannot be simple. \qedhere
	\end{enumerate}
\end{proof}

\section*{Acknowledgements}
The authors would like to thank the anonymous referees, whose comments greatly improved the presentation of this manuscript.
This work was partially supported by Universidad Nacional de San Luis, grants PI UNSL 03-1326, PROIPRO 03-2923, and PROINI 03-124, and Consejo Nacional de Investigaciones Cient\'ificas y T\'ecnicas grant PIP $11220220$ $100068$CO. 


\end{document}